\documentclass[11pt,leqno,a4paper]{article}

\usepackage{amsmath, amssymb, amsthm, amsfonts}
\usepackage{mathrsfs}
\usepackage{color}
\usepackage{hyperref}
\theoremstyle{plain}
\newtheorem{thm}{Theorem}
\newtheorem{cor}{Corollary}
\newtheorem{lem}{Lemma}

\theoremstyle{remark}
\newtheorem{rem}{Remark}

\renewcommand{\Re}{{\rm Re\,}}

\def\K{\mathop{\mbox{\bf\Large K}}}

\numberwithin{equation}{section}

\begin{document}

\title{The fastest possible continued fraction approximations of a class of functions}
\author{Xiaodong Cao*, Yoshio Tanigawa and Wenguang Zhai }
\date{}

\maketitle
\footnote[0]{* Corresponding author.}

\footnote[0]{\emph{E-mail address}: caoxiaodong@bipt.edu.cn~(X.D. Cao), tanigawa@math.nagoya-u.ac.jp~(Y. Yanigawa), zhaiwg@hotmail.com~(W.G. Zhai)}



\footnote[0]{This work is supported by
the National Natural Science Foundation of China (Grant No.11171344) and the Natural
Science Foundation of Beijing (Grant No.1112010).}

\footnote[0]{
Xiaodong Cao:
Department of Mathematics and Physics,
Beijing Institute of Petro-Chemical Technology,
Beijing, 102617, P. R. China \\
Yoshio Tanigawa:
Graduate School of Mathematics,
Nagoya University,
Nagoya, 464-8602, Japan\\
Wenguang Zhai:
Department of Mathematics,
China University of Mining and Technology,
Beijing 100083, P. R. China}

\begin{abstract} The goal of this paper is to formulate  a systematical method for constructing the fastest possible continued fraction approximations of a class of functions. The main tools are the multiple-correction method, the generalized Mortici's lemma and the
Mortici-transformation. As applications, we will present some sharp inequalities, and the continued fraction expansions associated to the volume of the unit ball. In addition, we obtain a new continued fraction expansion of Ramanujan for a ratio of the gamma functions, which is showed to be the fastest possible. Finally, three conjectures are proposed.
\end{abstract}
\section{Introduction}
Let $f(x)$ be a function defined on $(0,+\infty)$ to be approximated. We suppose that there exists a fixed positive integer $\nu$ and a constant $c\neq 0$ such that
\begin{align}
\lim_{x\rightarrow +\infty}x^{\nu}f(x)=c.\label{Condition-1}
\end{align}

In this case, we say that the function $f(x)$ is of order $x^{-\nu}$ when $x$ tends to infinity, and denote
\begin{align}
 \mathrm{R}(f(x)):=\nu,
\end{align}
where $\nu$ is the exponent of $x^{\nu}$. For convenience, $\mathrm{R}(0)$ is stipulated to be infinity. Hence, $\mathrm{R}(f(x))$ characterizes the rate of convergence for $f(x)$ as $x$ tends to infinity. From~\eqref{Condition-1},
there exists a large positive number $X_0$ such that $f(x)/c>0$
when $x>X_0$.

%Without loss of generality, we assume $c=1$.
\bigskip
In analysis, approximation theory, applied mathematics, etc., we often need to investigate the rational function approximation problem. Let $\frac{P_l(x)}{Q_m(x)}$ be an approximation to $f(x)$ as $x$ tends to infinity, where $P_l(x)$ and $Q_m(x)$ are polynomials in $x$. Quite similarly to the rational approximation problem for an irrational number, in order to find a better approximation to $f(x)$, we have to increase the degrees of both $P_l(x)$ and $Q_m(x)$. The main interest in this paper is to try to look for the fastest possible continued
fraction approximation or guess its approximation structure for $f(x)$ as $x$ tends to infinity.
\bigskip

The paper is organized as follows. In Sec.~2, we mainly introduce a definition to classify the continued fraction. In Sec.~3, we will prepare two preliminary lemmas for later use. In Sec.~4, we first develop further the previous
\emph{multiple-correction method}. Secondly, we introduce a transformation named as~\emph{Moritici-transformation} to change a kind of continued fraction approximation problem. In addition, we also give its \emph{Mathematica} program for the reader's convenience. Thirdly, similarly to Taylor's formula, we introduce two definitions of the \emph{formal Type-I and Type-II continued fraction approximation of order $k$} for a function, and the formal continued fraction expansion, respectively. This section constitutes the main part of this paper. To illustrate our method formulated in Sec.~4, in Sec.~5 we use the volume of the unit ball as an example to present some new inequalities. In Sec.~6, we test the well-known generalized Lord Brouncker's continued fraction formula, and show that it is the fastest possible. We also give some applications for the continued fraction formula involving the volume of the unit ball. In Sec.~7, we will use a continued fraction formula of Ramanujan to illustrate how to get the fastest possible  form of the continued fraction expression. In Sec.~8, we explain how to guess the fastest possible continued fraction expansions, and give three conjectures associated to the special rate of gamma functions. In the last section, we analyze the related perspective of research in this direction.

\section{Notation and definition}
Throughout the paper, we use the notation $\lfloor x\rfloor$ to denote the largest integer not exceeding $x$. The notation
$P_k(x)$(or $Q_k(x)$) means a polynomial of degree $k$
in $x$. We will use the $\Phi(k;x)$ to denote a polynomial of degree $k$ in $x$ with the leading coefficient equals one, which may be different at each occurrence. While, the notation $\Psi(k;x)$
means a polynomial of degree $k$ in $x$ with all coefficients non-negative, which may be different at each occurrence.
Let $(a_n)_{n\ge 1}$ and $(b_n)_{n\ge 0}$ be two sequences of real numbers with $a_n\neq 0 $ for all $n\in\mathbb{N}$ . The generalized continued fraction
\begin{align}
\tau=b_0+\frac{a_1}{b_1+\frac{a_2}{b_2+\ddots}}=b_0+
\begin{array}{ccccc}
a_1 && a_2 &       \\
\cline{1-1}\cline{3-3}\cline{5-5}
 b_1 & + & b_2 & + \cdots
\end{array}
=b_0+\K_{n=1}^{\infty}
\left(\frac{a_n}{b_n}\right)
\end{align}
is defined as the limit of the $n$th approximant
\begin{align}
\frac{A_n}{B_n}=b_0+\K_{k=1}^{n}\left(\frac{a_k}{b_k}\right)
\end{align}
as $n$ tends to infinity. The canonical numerators $A_n$ and denominators $B_n$ of the approximants satisfy the recurrence relations~(see [8, p. 105])
\begin{align}
A_{n+2}=b_{n+2}A_{n+1}+a_{n+2}A_n,\quad B_{n+2}=b_{n+2}B_{n+1}+a_{n+2}B_n
\end{align}
with the initial values $A_0=b_0, B_0=1, A_1=b_0b_1+a_1$ and
$B_1=b_1$.
\bigskip

To describe our method clearly, we will introduce two definitions as follows.
\bigskip

\noindent {\bf Definition 1.} Let $c_0\neq 0$, and $x$ be a free variable. Let  $(a_n)_{n=0}^{\infty}$, $(b_n)_{n=0}^{\infty}$ and $(c_n)_{n=0}^{\infty}$ be three real sequences. The formal continued fraction
\begin{align}
\frac{c_0}{\Phi(\nu;x)+\K_{n=0}^{\infty}\left(\frac{a_n}{x+b_n}\right)}
\end{align}
is said to be a \emph{Type-I} continued fraction. While,
\begin{align}
\frac{c_0}{\Phi(\nu;x)+\K_{n=0}^{\infty}\left(\frac{a_n}{x^2+b_n x+c_n}\right)}
\end{align}
is said to be a \emph{Type-II} continued fraction.
\bigskip

\begin{rem} The \emph{Type-I} and \emph{Type-II} are two kinds of fundamental structures we often meet. Certainly, we may define other-type continued fraction. Because of their complexity, in this paper we will not discuss the involved problems.
\end{rem}
\bigskip

\noindent {\bf Definition 2.} If the sequence $(b_n)_{n=0}^{\infty}$ is a constant sequence $(b)_{n=0}^{\infty}$ in the \emph{Type-I}~( or \emph{Type-II})
continued fraction, we call the number $\omega=b$~(or $\omega=\frac b2$) the $\mathrm{MC}$-point for the corresponding continued fraction.
We use $\hat{x}=x+\omega$ to denote the $\mathrm{MC}$-shift of $x$.
\bigskip

If there exists the $\mathrm{MC}$-point, we have the following \emph{simplified form}
\begin{align}
\frac{c_0}{\Phi_1(\nu;\hat{x})+\K_{n=0}^{\infty}\left(\frac{a_n}
{\hat{x}}\right)}\quad
\mbox{or}\quad
\frac{c_0}{\Phi_1(\nu;\hat{x})+\K_{n=0}^{\infty}\left(\frac{a_n}
{\hat{x}^2+d_n}\right)},\label{canonical-form}
\end{align}
where $d_n=c_n-\frac{b^2}{4}.$

\section{Two preliminary lemmas }
%A generalized Mortici's lemma
Mortici~\cite{Mor1} established a very useful tool for measuring the rate of convergence, which claims that a sequence $(x_n)_{n\ge 1}$ converging to zero is the fastest possible when the difference $(x_n-x_{n+1})_{n\ge 1}$ is the fastest possible. Since then, Mortici's lemma has been effectively applied in many papers such as~\cite{CXY,Cao1,CY,CW,Mor2,Mor3,Mor4,MCL}.
The following lemma is a generalization of Mortici's lemma. For details, readers may refer to~\cite{Cao2}.
\begin{lem} If $\lim_{x\rightarrow+\infty}f(x)=0$, and there exists the limit
\begin{align}
 \lim_{x\rightarrow+\infty}x^\lambda\left(f(x)-f(x+1)\right)=l\in
 \mathbb{R},
\end{align}
with $\lambda>1$, then
\begin{align}
 \lim_{x\rightarrow+\infty}x^{\lambda-1}f(x)=\frac{l}{\lambda-1}.
\end{align}
\end{lem}

\bigskip

In this paper, we will use the following simple inequality,
which is a consequence of Hermite-Hadamard inequality.

\begin{lem} Let $f$ be twice differentiable
with $f''$ continuous.  If $f''(x)>0$, then
\begin{align}
\int_{a}^{a+1}f(x) dx > f(a+1/2).\label{LEM3}
\end{align}
\end{lem}

\section{The multiple-correction, the Mortici-transformation and the formal continued fraction expansion}

\subsection{The multiple-correction method}
In this subsection, we will develop further the previous \emph{multiple-correction method} formulated in~\cite{Cao1,Cao2}. For some applications of this method, reader may refer to~\cite{CXY,Cao2,CY,CW}. In fact, the
\emph{multiple-correction method} is a recursive algorithm, and one of its advantages is that by repeating correction-process we always can accelerate the convergence. More precisely, every non-zero coefficient plays an important role in accelerating the convergence. The
\emph{multiple-correction method} consists of the following several steps.
\bigskip

\noindent {\bf(Step 1) The initial-correction.} The initial-correction is
vital. Determine the initial-correction $\Phi_0(\nu;x)$ such that
\begin{align}
\mathrm{R}\left(f(x)-\frac{c}{\Phi_0(\nu;x)}\right)=
\max_{\Phi(\nu;x)}\mathrm{R}
\left(f(x)-\frac{c}{\Phi(\nu;x)}\right).
\end{align}
\bigskip

\noindent {\bf(Step 2) The first-correction.} If there exists a real number $\kappa_0$ such that
\begin{align}
\mathrm{R}\left(f(x)-\frac{c}{\Phi_0(\nu;x)+\frac{\kappa_0}{x}}\right)
>\mathrm{R}\left(f(x)-\frac{c}{\Phi_0(\nu;x)}\right),
\end{align}
then we take the first-correction $\mathrm{MC}_1(x)=\frac{\kappa_0}{x+\lambda_0}$ with
\begin{align}
\lambda_0=\max_{\lambda}\mathrm{R}\left(f(x)-\frac{c}{\Phi_0(\nu;x)
+\frac{\kappa_0}{x+\lambda}}\right).
\end{align}
In this case, the first-correction has the form \emph{Type-I}. Otherwise, we take the first-correction $\mathrm{MC}_1(x)$ in the form \emph{Type-II}, i.e. $\mathrm{MC}_1(x)=\frac{\kappa_0}{x^2+\lambda_{0,1}x+\lambda_{0,2}}$ such that
\begin{align}
(\kappa_0,\lambda_{0,1},\lambda_{0,2})=
\max_{\kappa,\lambda_1,\lambda_2}\mathrm{R}
\left(f(x)-\frac{c}{\Phi_0(\nu;x)
+\frac{\kappa}{x^2+\lambda_1 x+\lambda_2}}\right).
\end{align}
If $\kappa_0=0$, we stop the correction-process, which means that the rate of convergence can not be further improved only by making use of \emph{Type-I} or \emph{Type-II} continued fraction structure.
\bigskip

\noindent {\bf(Step 3) The second-correction to the $k$th-correction.} If $\mathrm{MC}_1(x)$ has the form \emph{Type-I}, we take the second-correction
\begin{align}
\mathrm{MC}_2(x)=\frac{\kappa_0}
{x+\lambda_0+\frac{\kappa_1}{x+\lambda_1}},
\end{align}
which satisfies
\begin{align}
(\kappa_1,\lambda_1)=\max_{\kappa,\lambda}\mathrm{R}
\left(f(x)-\frac{c}{\Phi_0(\nu;x)
+\frac{\kappa_0}
{x+\lambda_0+\frac{\kappa}{x+\lambda}}}\right).
\end{align}
Similarly to the first-correction, if $\kappa_1=0$, we stop the correction-process.

If $\mathrm{MC}_1(x)$ has the form \emph{Type-II}, we take the second-correction
\begin{align}
\mathrm{MC}_2(x)=\frac{\kappa_0}
{x^2+\lambda_{0,1}x+\lambda_{0,2}+\frac{\kappa_1}
{x^2+\lambda_{1,1}x+\lambda_{1,2}}},
\end{align}
such that
\begin{align}
(\kappa_1,\lambda_{1,1},\lambda_{1,2})=
\max_{\kappa,\lambda_{1},\lambda_{2}}\mathrm{R}
\left(f(x)-\frac{c}{\Phi_0(\nu;x)
+\frac{\kappa_0}
{x^2+\lambda_{0,1}x+\lambda_{0,2}+\frac{\kappa}
{x^2+\lambda_{1}x+\lambda_{2}}}}\right).
\end{align}
If $\kappa_1=0$, we also need to stop the correction-process.

If we can continue the above correction-process to determine the $k$th-correction function $\mathrm{MC}_k(x)$ until some $k^*$ you want, then one may use a recurrence relation to determine $\mathrm{MC}_k(x)$. More precisely, in the case of \emph{Type-I} we choose
\begin{align}
\mathrm{MC}_k(x)=\K_{j=0}^{k-1}
\left(\frac{\kappa_j}{x+\lambda_j}\right)
\end{align}
such that
\begin{align}
(\kappa_{k-1},\lambda_{k-1})=\max_{\kappa,\lambda}\mathrm{R}
\left(f(x)-\left(\begin{array}{ccccccc}
c& &\kappa_0 & & \kappa_{k-2} &      & \kappa \\
\cline{1-1}\cline{3-3}\cline{5-5}\cline{7-7}
\Phi_0(\nu;x) &+&x+\lambda_0 & +\cdots+ & x+\lambda_{k-2} & + &x+\lambda
\end{array}\right)\right).
\end{align}
While, in the case of \emph{Type-II} we take
\begin{align}
\mathrm{MC}_k(x)=\K_{j=0}^{k-1}\left(\frac{\kappa_j}
{x^2+\lambda_{j,1}x+\lambda_{j,2}}\right),
\end{align}
which satisfies
\begin{align}
(\kappa_{k-1},\lambda_{k-1,1},\lambda_{k-1,2})
=\max_{\kappa,\lambda_{1},\lambda_{2}}\mathrm{R}
\left(f(x)-G(\kappa,\lambda_1,\lambda_2;x)\right),
\end{align}
where
\begin{align*}
G(\kappa,\lambda_1,\lambda_2;x):=\begin{array}{ccccccc}
c& &\kappa_0 & & \kappa_{k-2} &      & \kappa \\
\cline{1-1}\cline{3-3}\cline{5-5}\cline{7-7}
\Phi_0(\nu;x) &+& x^2+\lambda_{0,1}x+\lambda_{0,2} & +\cdots+ & x^2+\lambda_{k-2,1}x+\lambda_{k-2,2} & + &x^2+\lambda_1 x+\lambda_2
\end{array}.
\end{align*}
\bigskip

Note that in the case of both \emph{Type-I} and \emph{Type-II} continued fraction approximation, if $\kappa_{k-1}=0$, we must stop the correction-process. In other words, to improve the rate of convergence, we need to choose some more complex continued fraction structure instead of it.
\bigskip

\begin{rem}
Sometimes, we need to consider its equivalent forms. For example,
the Stirling's formula reads (See, e.g. [1, p. 253])
\begin{align}
  \Gamma(x+1)\sim \sqrt{2\pi x} \left(\frac xe\right)^{x},\quad x\rightarrow +\infty,\label{Stirling's formula}
\end{align}
which is equivalent to
\begin{align}
\lim_{x\rightarrow \infty}x^3 f(x)=1,\label{Ramanujan-Type}
\end{align}
where
\begin{align}
f(x)=8 \pi^3\left(\frac xe\right)^{6x}\Gamma^{-6}(x+1).
\end{align}
From the above asymptotic formula, we may study Ramanujan-type continued fraction approximation for the gamma function. For more details, see Cao~\cite{Cao2} or next section. Moreover, we note that~\eqref{Stirling's formula} has many equivalent forms. Hence, it is not difficult to see that the equivalent transformation of a practical problem influences directly the initial-correction and final continued fraction approximation.
\end{rem}

\begin{rem}
If $\nu$ is a negative integer, our method is still efficient, i.e. we may consider the reciprocal of $f(x)$.
\end{rem}
\begin{rem}
For comparison, we use the mathematical notation ``$\mathrm{R}$" and ``$\max$" in the above definition, which make the method more clearly.
\end{rem}
\subsection{The Mortici-transformation}
In this subsection we will explain how to look for all the related coefficients in $\Phi_0(\nu;x)$ and $\mathrm{MC}_k(x)$. If we can expand $f(x)$ into a power series in terms of $1/x$ easily, then it is not difficult to determine $\Phi_0(\nu;x)$ and $\mathrm{MC}_k(x)$. Similarly, if we may expand the difference $f(x)-f(x+1)$ into a power series in terms of $1/x$, by the generalized Moritici's lemma we also can find $\Phi_0(\nu;x)$ and $\mathrm{MC}_k(x)$, e.g. the Euler-Mascheroni
constant, the constants of Landau, the constants of Lebesgue, etc.~(See~\cite{Cao1}). However, in many cases the previous two approaches are not very efficient, e.g. gamma function~(see, Remark 2) and  the ratio of the gamma functions~(for example, see Sec.~7 below). Instead, we may employ the following method to achieve it.
\bigskip

First, we introduce the $k$th-correction relative error sequence $(E_k(x))_{k\ge 0}$ as follows
\begin{align}
&f(x)=\frac{c}{\Phi_0(\nu;x)}
\exp\left(E_0(x)\right),\label{E0-def}\\
&f(x)=\frac{c}{\Phi_0(\nu;x)+\mathrm{MC}_k(x)}
\exp\left(E_k(x)\right),\quad k\ge 1,\label{Ek-def}
\end{align}
where $\Phi_0(k;x)$ is a polynomial of degree $\nu$ in $x$ with the leading coefficient equals one, to be specified below.

It is easy to verify that
\begin{align*}
&f(x)-\frac{c}{\Phi_0(\nu;x)}=\frac{c}{\Phi_0(\nu;x)}\left(
\exp\left(E_0(x)\right)-1
\right),\\
&f(x)-\frac{c}{\Phi_0(\nu;x)+\mathrm{MC}_k(x)}=
\frac{c}{\Phi_0(\nu;x)+\mathrm{MC}_k(x)}
\left(
\exp\left(E_k(x)\right)-1
\right),\quad k\ge 1.
\end{align*}

It is well-known that
\begin{align*}
\lim_{t\rightarrow 0}\frac{\exp(t)-1}{t}=1,
\end{align*}
by $\lim_{x\rightarrow\infty}E_k(x)=0$ we obtain
\begin{align}
&\mathrm{R}\left(f(x)-\frac{c}{\Phi_0(\nu;x)}\right)
=\nu+\mathrm{R}\left(E_0(x)\right),\\
&\mathrm{R}\left(f(x)-\frac{c}{\Phi_0(\nu;x)+\mathrm{MC}_k(x)}\right)
=\nu+\mathrm{R}\left(E_k(x)\right),\quad k\ge 1.
\end{align}
In this way, we turn the problem to solve
$\mathrm{R}\left(E_k(x)\right)$.%reduce

Take the logarithm of~\eqref{E0-def} and~\eqref{Ek-def}, respectively, we deduce that
\begin{align*}
&\ln \frac{f(x)}{c}=-\ln\left(\Phi_0(\nu;x)\right)+E_0(x),\\
&\ln \frac{f(x)}{c} =-\ln\left(\Phi_0(\nu;x)+\mathrm{MC}_k(x)\right)+E_k(x),\quad k\ge 1.
\end{align*}

Next, let us consider the difference
\begin{align}
E_0(x)-E_0(x+1)=&\ln\frac{f(x)}{f(x+1)}
+\ln\frac{\Phi_0(\nu;x)}{\Phi_0(\nu;x+1)},\\
E_k(x)-E_k(x+1)=&\ln\frac{f(x)}{f(x+1)}
+\ln\frac{\Phi_0(\nu;x)+\mathrm{MC}_k(x)}{\Phi_0(\nu;x+1)
+\mathrm{MC}_k(x+1)},\quad k\ge 1.
\end{align}
By Lemma 1~(the generalized Moritici's lemma), we have
\begin{align}
\mathrm{R}\left(E_k(x)\right)=\mathrm{R}\left(E_k(x)-E_k(x+1)\right)-1.
\end{align}

Finally, if set $\mathrm{MC}_0(x)\equiv 0$, then we attain the following useful tool.
\bigskip

\begin{lem} Let $f(x)$ satisfy~\eqref{Condition-1}. Under the above notation, we have
\begin{align}
&\mathrm{R}\left(f(x)-\frac{c}{\Phi_0(\nu;x)+\mathrm{MC}_k(x)}\right)
\label{Mortici-transformation}\\
=&\nu-1+\mathrm{R}\left(\ln\frac{f(x)}{f(x+1)}
+\ln\frac{\Phi_0(\nu;x)+\mathrm{MC}_k(x)}{\Phi_0(\nu;x+1)
+\mathrm{MC}_k(x+1)}\right),\quad k\ge 0.\nonumber
\end{align}
\end{lem}
The idea of Lemma 3 is first originated from Mortici~\cite{Mor1}, which will be called a \emph{Mortici-transformation}. We would like to stress that \emph{Mortici-transformation} implies the following assertion
\begin{align}
&\max_{\kappa,\lambda \ (\mbox{\emph{\footnotesize or}}~ \kappa,\lambda_1,\lambda_2)}\mathrm{R}\left(f(x)-\frac{c}
{\Phi_0(\nu;x)+\mathrm{MC}_k(x)}\right)
\label{Mortici-transformation-1}\\
=&\max_{\kappa,\lambda \ (\mbox{\emph{\footnotesize or}}~ \kappa,\lambda_1,\lambda_2)}\mathrm{R}\left(\ln\frac{f(x)}{f(x+1)}
+\ln\frac{\Phi_0(\nu;x)+\mathrm{MC}_k(x)}{\Phi_0(\nu;x+1)
+\mathrm{MC}_k(x+1)}\right),\quad k\ge 0.\nonumber
\end{align}
In the sequel, we will use this relation many times. For the sake of simplicity, we will always assume that the difference
\begin{align}
\ln\frac{f(1/z)}{c}- \ln\frac{f(1/z +1)}{c} =\ln\frac{f(1/z)}{f(1/z+1)}\label{Condition-2}
\end{align}
is an analytic function in a neighborhood of point $z=0$.
\bigskip

For the reader's convenience, we would like to give the complete \emph{Mathematica} program for finding
all the coefficients in $\Phi_0(\nu;x)$ and $\mathrm{MC}_k(x)$ by making use of~\emph{Mortici-transformation}.
\bigskip

{\bf (i)}. First, let the function $MT[x]$ be defined by
\begin{align*}
MT[x]:=\ln\frac{f(x)}{f(x+1)}
+\ln\frac{\Phi_0(\nu;x)+\mathrm{MC}_k(x)}{\Phi_0(\nu;x+1)
+\mathrm{MC}_k(x+1)}.
\end{align*}

{\bf (ii)}. Then we manipulate the following \emph{Mathematica} command to
expand $MT[x]$ into a power series in terms of $1/x$:
\begin{align}
\text{Normal}[\text{Series}
[MT[x]
\text{/.}~ x\rightarrow
1/u, \{u,0,l_k\}]]\text{/.}~ u\rightarrow 1/x~(\text{// Simplify})
\label{MT-Mathematica-Program}
\end{align}
We remark that the variable $l_k$ needs to be suitable chosen according to the different function.
\bigskip

{\bf (iii)}. Taking out the first some coefficients in the above power series, then we enforce them to be zero, and finally solve the related coefficients successively.

\begin{rem}
Actually, once we have found $\mathrm{MC}_k(x)$, \eqref{MT-Mathematica-Program} can be used again to determine the rate of convergence. In addition, we can apply it to check the general term formula for $\mathrm{MC}_k(x)$.
\end{rem}

\subsection{The formal continued fraction expansion}

Similarly to Taylor's formula, if the $k$th-correction $\mathrm{MC}_k(x)$ for $f(x)$ has the \emph{Type-I (or the Type-II)} structure, then we may construct the \emph{formal Type-I (or Type-II) continued fraction approximation of order $k$} for $f(x)$ as follows:
\begin{align}
CF_k(f(x)):= \frac{1}{\Phi_0(\nu;x)+\mathrm{MC}_k(x)},\quad k\ge 0.
\end{align}

For example, Euler-Mascheroni
constant has the formal \emph{Type-I} continued fraction approximation of order $k$, while both Landau's constants and Lebesgue's constants have the formal \emph{Type-II} continued fraction approximation of order $k$. For details, readers may refer to~\cite{Cao1}.
\bigskip

{\bf Example 1.} Let $f(x)=\frac{\Gamma^4(x+\frac 14)}{\Gamma^4(x+1)}$. Then $CF_k(f(x))$ is the \emph{Type-I}, its $\mathrm{MC}$-point $\omega$ equals $\frac 18$~(i.e. $\lambda_m\equiv \frac 18$), and
\begin{align}
CF_k(f(x))=\frac{1}{(x+\frac 18)^3+\frac {7}{128} (x+\frac 18)+\mathrm{MC}_k(x)},
\end{align}
where $
(\kappa_0,\kappa_1,\kappa_2,\ldots)=\left(-\frac{189}{32768}, \frac{1483}{2688}, \frac{8923253}{15945216}, \frac{
10136617131375}{6775390309888}, \frac{2313439127848201}{1428763287038592},\ldots\right).
$

{\bf Example 2.} Let $G_{\eta}(x)$ be defined by \eqref{G-eta-Def}
below. In the case of $\eta=\frac 12$, $CF_k\left(G_{\eta}^2(x)\right)$ is the \emph{Type-II}, for details see Corollary 2 in Sec.~7. If $\eta\neq \frac 12$, then $CF_k\left(G_{\eta}^2(x)\right)$ is the \emph{Type-I}, and it has not $\mathrm{MC}$-point. We have
\begin{align}
CF_2\left(G_{\eta}^2(x)\right)=\frac{1}{x^2+2\eta(1-\eta)x+2 \eta^2(\eta-1)^2 +\mathrm{MC}_2(x)},
\end{align}
where $\kappa_0=-\frac 13\eta^2 (1-\eta)^2(2 \eta-1)^2$, $\lambda_0=\frac{(2 \eta - 1)^2}{8} + \frac 14 - \frac{3}{8 (2 \eta-1)^2}$, $\kappa_1=\frac{1}{64} \left( (2 \eta - 3)^2(2 \eta + 1)^2  + 10 + \frac{45}{(2 \eta-1)^4}\right)$, and $\lambda_1=\frac{(2 \eta - 1)^2}{24} + \frac 14 + \frac{3}{8 (2 \eta-1)^2} -\frac{( \eta-2)^2 (\eta+1)^2 (2 \eta- 1)^2}{6(2 \eta-1)^4\kappa_1}$.
\bigskip

If we rewrite $CF_k(f(x))$ in a rational function of the form $\frac{P_r(x)}{Q_s(x)}$, then $s=k+\nu$ in the case of  \emph{Type-I}, and $s=2k+\nu$ in the case of \emph{Type-II}. If we let $\mathrm{R}\left(f(x)-CF_k(f(x))\right)=K$, then
\begin{align}
f(x)=CF_k(f(x))+O\left(x^{-K}\right),\quad x\rightarrow\infty.
\end{align}
\bigskip

Let $\theta_0=0$ or $1$. A lot of computations reveal that if $CF_k(f(x))$
is the \emph{Type-I}, then $K=2k+2\nu+1+\theta_0$, and
$K=4k+2\nu+1+\theta_0$ in the case of \emph{Type-II}, respectively.
\bigskip

For a suitable ``not very large" positive integer $k$, by using of \emph{Mortici-transfomation} and~\eqref{MT-Mathematica-Program}, we may get the rate of convergence for $f(x)-CF_k(f(x))$ when $x$ tends to infinity. Moreover, by making use of telescoping method, Hermite-Hadamard inequality, etc, sometimes we can prove sharp double inequalities of $f(x)-CF_k(f(x))$ for as smaller $x$ as possible. We will give an example in Sec.~5.
\bigskip

Now let $k$ tend to $\infty$, we get the \emph{formal Type-I (or Type-II) continued fraction expansion} for $f(x)$, or shortly write \begin{align}
f(x)\sim CF(f(x)):=CF_{\infty}(f(x)),\quad x\rightarrow \infty.
\end{align}

In some cases, we can test and guess further the general term of
$CF(f(x))$. Here we need to apply some tools in number theory,
difference equation, etc. We will show some examples in Sec.~7.
% occasionally.
\bigskip

For the formal continued fraction expansion, we are often concerned with the following two main problems.
\bigskip

\noindent {\bf Problem 1.}  Determine the domains of convergence for the formal continued fraction expansion $CF(f(x))$. We may refer to two very nice books:~L. Lorentzen and H. Waadeland\cite{LW}, and A. Cuyt, V.B. Petersen, B. Verdonk, H. Waadeland, W.B. Jones~\cite{CPV}, or some other classical books cited in there.
\bigskip

\noindent {\bf Problem 2.} Prove an identity for as the large domains as possible. That is, based on Problem 1, to determine the intervals $\mathrm{I}$ such that $f(x)= CF(f(x))$ for all $x\in \mathrm{I}$. For example, with the help of continued fraction theory, hypergeometric series, etc., we hope at least to find a interval $(x_0,\infty)\subset \mathrm{I}$ for some $x_0>0$. Certainly, we may extend it to a complex domain. However, in this paper we will not investigate this topic.
\bigskip

On one hand, to determine all the related coefficients, we often use an appropriate symbolic computation software, which needs a huge of computations. On the other hand, the exact expressions at each occurrence also takes a lot of space. Hence, in this paper we omit some related details for space limitation.

\begin{rem}
From the above discussion, we observe that for a specific function, except a huge of computations, probably only such two kinds of structures can not provide ``good continued fraction approximation''. In addition, in the theory of classical continued fraction, even if there is a continued fraction expansion for a given function, we often do not know whether it is the fastest possible or best possible. Generally speaking, for a given continued fraction, finding the rate of convergence for the $k$th approximant is not always easy.
\end{rem}
\section{The volume of the unit ball}
It is well-known that the volume of the unit ball in $\mathbb{R}^n$ is
\begin{align}
\Omega_n=\frac{\pi^{\frac{n}{2}}}{\Gamma(\frac n2+1)}.\label{Omega-n-Def}
\end{align}
Many authors have investigated the inequalities about the $\Omega_n$, e.g. see~\cite{Al1,Al2,Al3,AVV,AQ,BH,Gao,KR,Me,Mor3,Mor4,QV,Sch} and references therein.

Chen and Li~\cite{CL} proved~($a=\frac e2,~b=\frac 13$):
\begin{align}
\frac{1}{\sqrt{\pi(n+a)}}\left(\frac{2\pi e}{n}\right)^{\frac n2}\le \Omega_n<\frac{1}{\sqrt{\pi(n+b)}}\left(\frac{2\pi e}{n}\right)^{\frac n2}.
\end{align}

Recently, Mortici~[25, Theorem 3] showed that for every integer $n\ge 3$ in the left-hand side and $n\ge 1$ in the right-hand side, then we have the following Gosper-type inequalities:
\begin{align}
\frac{1}{\sqrt{\pi(n+\theta(n))}}\left(\frac{2\pi e}{n}\right)^{\frac n2}\le \Omega_n<\frac{1}{\sqrt{\pi(n+\vartheta(n))}}\left(\frac{2\pi e}{n}\right)^{\frac n2},
\end{align}
where
\begin{align*}
\theta(n)=\frac 13+\frac{1}{18 n}-\frac{31}{810 n^2},\quad \vartheta(n)=\theta(n)-\frac{139}{9720 n^3}.
\end{align*}
Now we let
\begin{align}
V(x)=\frac{\pi^x}{\Gamma(x+1)},\quad x>0.\label{V-Def}
\end{align}
Let us imagine that if $\frac{1}{\Gamma(x+1)}\sim H(x)$ when $x$ tends to infinity, then $V(x)$ has an asymptotic formula of the form
$\pi^x H(x)$. In this sense, by Remark 2 and Remark 3, it suffices to consider the asymptotic formula for the gamma function. In fact, we note that both $f(x)$ and $1/f(x)$ have the same $k$th-correction $MC_k(x)$.

From~\eqref{Ramanujan-Type}, we introduce the relative error sequence $(E_k(x))_{k\ge 0}$ to be defined by
\begin{align}
f(x):=&8 \pi^3\left(\frac xe\right)^{6x}\Gamma^{-6}(x+1)=\frac{\exp(E_0(x))}{\Phi_0(x)},
\label{E0-Gamma-Def}\\
f(x):=&\frac{\exp(E_k(x))}{\Phi_0(x)+\mathrm{MC}_k(x)},\quad k\ge 1,\label{Ek-Gamma-Def}
\end{align}
where $\Phi_0(x)= x^3+\frac 12 x^2+\frac 18 x+\frac {1}{240}$, and
\begin{align}
\mathrm{MC}_k(x)=\K_{j=0}^{k-1}\left(\frac{\kappa_j}{x+\lambda_j}
\right),
\end{align}
here $\kappa_0=-\frac{11}{1920}, \lambda_0=\frac{79}{154},
\kappa_1=\frac{459733}{711480}, \lambda_1=-\frac{1455925}{70798882},\ldots$. We stress that $\Phi_0(x)$ was claimed  first by Ramanujan~\cite{Ram}, and some more coefficients may be founded in~\cite{Cao2}. By employing Lemma 1, \eqref{Ek-difference-Exp}~(see below) and~\eqref{MT-Mathematica-Program}, it is not difficult to verify that
\begin{align}
&\lim_{x\rightarrow\infty}x^{4}E_0(x)=\frac{11}{1920}:=C_0,\\
&\lim_{x\rightarrow\infty}x^{6}E_1(x)=-\frac{459733}{124185600}:=C_1.
\end{align}
The following theorem tells us how to improve the above results and obtain some sharper estimates for $E_0(x)$ and $E_1(x)$.

\begin{thm} Let $E_0(x)$ and $E_1(x)$ be defined as~\eqref{E0-Gamma-Def} and~\eqref{Ek-Gamma-Def}, respectively.

(i) For every real number $x\ge 6$ in the left-hand side and $x\ge 12$ in the right-hand side, we have
\begin{align}
\frac{11}{1920}\frac{1}{(x+3)^4}<E_0(x)<\frac{11}{1920 (x-5)^4}.
\end{align}

(ii) For every real number $x\ge 9$ in the left-hand side and $x\ge 10$ in the right-hand side, then
\begin{align}
-\frac{459733}{124185600}\frac{1}{(x-2)^6}<E_1(x)
<-\frac{459733}{124185600}\frac{1}{(x+2)^6}.
\end{align}
\end{thm}
\proof We use the idea of Theorem 2 in~\cite{XY} or Theorem 1 in~\cite{CXY}. Let $G_k(x)=E_k(x)-E_k(x+1)$ for $k\ge 0$. We will employ the telescoping method. It follows from $\lim_{x\rightarrow\infty}E_k(x)=0$ that
\begin{align}
E_k(x)=\sum_{m=0}^{\infty}G_k(x+m), \quad(k=0,1).\label{Telecsoping Identity}
\end{align}

If $g(\infty)=g'(\infty)=0$, it is not difficult to prove that
\begin{align}
g(x)=-\int_{x}^{\infty}g'(s) ds=
\int_{x}^{\infty}\left(\int_{s}^{\infty}g''(t) dt\right) ds.
\label{Integral-Change-Form}
\end{align}
Note that the convenience $\mathrm{MC}_0(x)=0$. By~\eqref{E0-Gamma-Def} and~\eqref{Ek-Gamma-Def}, we have
\begin{align}
&E_k(x)=-6\ln \Gamma(x+1)+2\ln 2\pi +6 x(\ln x-1)+ \ln(\Phi_0(x)+\mathrm{MC}_k(x)),\\
&G_k(x)=E_k(x)-E_k(x+1)=6\left(1-x \ln(1+\frac 1x)\right)+\ln\frac{\Phi_0(x)+\mathrm{MC}_k(x)}
{\Phi_0(x+1)+\mathrm{MC}_k(x+1)}.\label{Ek-difference-Exp}
\end{align}
By using \emph{Mathematica} software, we can check that if $x>0$, then
\begin{align}
&G_0''(x)- \frac{11}{16 x^7}+ \frac{29}{8 x^8}=\frac{\Psi_1(15;x)}{96 x^8 (1 + x)^2 \Psi_2(12;x)}>0,\\
&G_0''(x)- \frac{11}{16 x^7}+ \frac{29}{8 x^8}-\frac{9031}{800 x^9}
=-\frac{\Psi_3(13;x)}{800 x^9 (1 + x)^2\Psi_2(12;x)}<0.
\end{align}
By~\eqref{Integral-Change-Form}, we get that when $x>0$,
\begin{align}
\frac{11}{480 x^5}- \frac{29}{336 x^6}<G_0(x)<\frac{11}{480 x^5}- \frac{29}{336 x^6}+\frac{9031}{44800 x^7}.\label{G0-bounds}
\end{align}
Similarly, if $x\ge \frac {1}{16}$, we have
\begin{align}
&G_1''(x)+ \frac{459733}{369600 x^9}- \frac{39872247}{
4743200 x^{10}}\\
=&-\frac{\Psi_4(20;x)(x-\frac{1}{16})+\frac{4490\dots 0225}{1441\cdots 5872}}{28459200 x^{10}\Psi_5(6;x)\left(\Psi_6(3;x)(x - \frac{1}{23})+\frac{7670381}{279841}\right)^2\Psi_7(8;x) }<0,\nonumber\\
&G_1''(x)+ \frac{459733}{369600 x^9}- \frac{39872247}{
4743200 x^{10}}+\frac{1092949825573}{32724285440 x^{11}}\\
=&\frac{\Psi_8(20;x)(x-\frac{1}{16})+\frac{2388\cdots 7275}{
1125\cdots 2624}}{490864281600 x^{11}\Psi_5(6;x)\left(\Psi_6(3;x)(x - \frac{1}{23})+\frac{7670381}{279841}\right)^2\Psi_7(8;x) }>0,\nonumber
\end{align}
and
\begin{align}
-\frac{459733}{20697600 x^7}+ \frac{13290749}{
113836800 x^8}-\frac{1092949825573}{2945185689600 x^9}<G_1(x)<-\frac{459733}{20697600 x^7}+ \frac{13290749}{
113836800 x^8}.\label{G1-bounds}
\end{align}
Now, combining \eqref{Telecsoping Identity}, \eqref{G0-bounds} and \eqref{G1-bounds}, we attain that
\begin{align}
&0<E_0(x)-\frac{11}{480}\sum_{m=0}^{\infty}\frac{1}{(x+m)^5}
+\frac{29}{336 } \sum_{m=0}^{\infty}\frac{1}{(x+m)^6}<\frac{9031}{44800 }\sum_{m=0}^{\infty}\frac{1}{(x+m)^7},\quad (x>0),\\
&-\frac{1092949825573}
{2945185689600 }\sum_{m=0}^{\infty}\frac{1}{(x+m)^9}<\\
&E_1(x)+\frac{459733}{20697600 }\sum_{m=0}^{\infty}\frac{1}{(x+m)^7}- \frac{13290749}{
113836800 }\sum_{m=0}^{\infty}\frac{1}{(x+m)^8}
<0,\quad (x>\frac{1}{16}).\nonumber
\end{align}

Let $j\ge 2$ and $x>\frac 12$. By Lemma 2, we obtain
\begin{align}
&\frac{1}{(j-1)x^{j-1}}=\int_x^{\infty}\frac{dt}{t^j}
<\sum_{m=0}^{\infty}\frac{1}{(x+m)^j}\\
&<\sum_{m=0}^{\infty}
\int_{x+m-\frac 12}^{x+m-\frac 12}\frac{dt}{t^j}
=\int_{x-\frac 12}^{\infty}\frac{dt}{t^j}=\frac{1}{(j-1)(x-\frac 12)^{j-1}}.\nonumber
\end{align}
By applying (5.22) and (5.24), under the condition $x\ge 6$ we have
\begin{align}
E_0(x)>&\frac{11}{480}\frac{1}{4x^4}-\frac{29}{336}\frac{1}{5 (x-\frac 12)^5}\\
=&\frac{11}{1920}\frac{1}{(x+3)^4}+
\frac{\Psi_1(7;x)(x-6)+2164192911}{13440 x^4 (3 + x)^4 (-1 + 2 x)^5}\nonumber\\
>&\frac{11}{1920}\frac{1}{(x+3)^4}.\nonumber
\end{align}
Similarly to (5.25), if $x\ge 12$, then
\begin{align}
E_0(x)<&\frac{11}{480}\frac{1}{4(x-\frac 12)^4}-\frac{29}{336}\frac{1}{5 x^5}+
\frac{9031}{44800}\frac{1}{6 (x-\frac 12)^6}\\
=&\frac{11}{1920 (x-5)^4}-\frac{\Psi_2(9;x)(x-12)+12561000435989768}{67200 (-1 + x)^4 x^5 (-1 + 2 x)^6}\nonumber\\
<&\frac{11}{1920 (x-5)^4}.\nonumber
\end{align}
This completes the proof of assertion (i). Finally, it is not difficult to check that if $x\ge 9$, then
\begin{align}
E_1(x)>&-\frac{459733}{20697600}\frac{1}{6 (x - 1/2)^6}+\frac{ 13290749}{
113836800}\frac{1}{7 x^7}  -\frac{1092949825573}{2945185689600}\frac{1}{8 (x - 1/2)^8}\\
=&-\frac{459733}{124185600}\frac{1}{(x - 2)^6}+\frac{\Psi_3(13;x)(x-9)+67733478135399363858702201}{736296422400 (-2 + x)^6 x^7 (-1 + 2 x)^8}\nonumber\\
>&-\frac{459733}{124185600}\frac{1}{(x - 2)^6},\nonumber
\end{align}
and if $x\ge 10$, we have
\begin{align}
E_1(x)<&-\frac{459733}{20697600}\frac{1}{6 x^6}+\frac{ 13290749}{
113836800}\frac{1}{7 (x - 1/2)^7}\\
=&-\frac{459733}{20697600}\frac{1}{6 (x+2)^6}-\frac{\Psi_4(11;x)(x-10)+470994290293217661904}{2390572800 x^6 (2 + x)^6 (-1 + 2 x)^7}\nonumber\\
<&-\frac{459733}{124185600}\frac{1}{(x + 2)^6},\nonumber
\end{align}
This will finish the proof of Theorem 1.\qed
\begin{thm} Assume $n\ge 24$, we have
the following Ramanujan-type inequalities
\begin{align}
\frac{1}{\sqrt \pi}\left(\frac{2\pi e}{n}\right)^{\frac n2}\frac{1}{\sqrt[6]{n^3+ n^2+ \frac n2+\frac{1}{30}}}<\Omega_n<\frac{1}{\sqrt \pi}\left(\frac{2\pi e}{n}\right)^{\frac n2}\frac{\exp\left(\frac{11}{720 (n-10)^4} \right)}{\sqrt[6]{n^3+ n^2+ \frac n2+\frac{1}{30}}}.
\end{align}
If $n\ge 20$, then
\begin{align}
&\frac{1}{\sqrt \pi}\left(\frac{2\pi e}{n}\right)^{\frac n2}\frac{1-\frac{459733}{11642400 (-4 + n)^6}}{\sqrt[6]{n^3+ n^2+ \frac n2+\frac{1}{30}-\frac{847}{9240 n + 9480}}}<\Omega_n\\
&<\frac{1}{\sqrt \pi}\left(\frac{2\pi e}{n}\right)^{\frac n2}\frac{1}{\sqrt[6]{n^3+ n^2+ \frac n2+\frac{1}{30}-\frac{847}{9240 n + 9480}}}.\nonumber
\end{align}
\end{thm}
\proof It follows from~\eqref{Omega-n-Def}, \eqref{E0-Gamma-Def} and \eqref{Ek-Gamma-Def} that
\begin{align}
\Omega_n=\frac{1}{\sqrt \pi}\left(\frac{2\pi e}{n}\right)^{\frac n2}\frac{\exp\left(\frac 16 E_0(\frac n2)\right)}{\sqrt[6]{8\Phi_0(\frac n2)}},\quad \Omega_n=\frac{1}{\sqrt \pi}\left(\frac{2\pi e}{n}\right)^{\frac n2}\frac{\exp\left(\frac 16 E_1(\frac n2)\right)}{\sqrt[6]{8\Phi_0(\frac n2)+8\mathrm{MC}_1(\frac n2)}}.
\end{align}
Now (5.29) follows from (5.10) and (5.31).

We begin to prove (5.30). It is well-known that $\exp(t)\ge 1+t$. When $n\ge 20$, by the inequality of the right-hand side in (5.11), we have the following trivial estimate
\begin{align}
\exp\left(\frac 16 E_1(\frac n2)\right)<1.
\end{align}
In addition, by the lower bound in (5.11), we get
\begin{align}
\exp\left(\frac 16 E_1(\frac n2)\right)>1+\frac 16 E_1(\frac n2)
>1-\frac{459733}{11642400 (-4 + n)^6},
\quad (n\ge 20).
\end{align}
Combining (5.31), (5.32) and (5.33) completes the proof of (5.30).\qed

Following the same approach as Theorem 2, it is not difficult to prove the following Ramanujan-type inequalities for the gamma function.
\begin{cor} Let $x\ge 12$. Then
\begin{align*}
\sqrt{\pi}\left(\frac x e\right)^x
\left(8x^3+4x^2+x+\frac{1}{30}\right)^{\frac 16}\exp\left(-\frac{11}{11520 (x-5)^4} \right)<\Gamma(x+1)<\sqrt{\pi}\left(\frac x e\right)^x
\left(8x^3+4x^2+x+\frac{1}{30}\right)^{\frac 16}.
\end{align*}
\end{cor}
\begin{rem}
It should is noted that the method described in Theorem 1 and 2 also can be used to look for $CF_k(F(x))$, and prove some inequalities involving the ratio of gamma functions.
\end{rem}
\begin{rem}
We will give some other results involving $\Omega_n$ in the subsection 6.2.
\end{rem}

\section{Lord Brouncker's continued fraction formula}
\subsection{Lord Brouncker's continued fraction formula}
The following formula is taken from Corollary 1 of Berndt~[8, p. 145], which was first proved by Bauer~\cite{Bau} in 1872.
\begin{lem} If $\Re x>0$, then
\begin{align}
\frac{\Gamma^2(\frac 14(x+1))}{\Gamma^2(\frac 14(x+3))}
=\frac{4}{x+}\K_{m=1}^{\infty}\left(\frac{(2m-1)^2}{2x}\right)
.\label{General
-Lord Brouncker}
\end{align}
\end{lem}
By taking $x=4n+1$ in the above formula, we obtain the so-called Lord Brouncker's continued fraction formula
\begin{align}
q(n):=\frac{\Gamma^2(n+\frac 12)}{\Gamma^2(n+1)}=
\frac{4}{4n+1+\frac{1^2}{2(4n+1)+\frac{3^2}{2(4n+1)
+\frac{5^2}{2(4n+1)+\ddots}}}}.\label{Brouncker's formula}
\end{align}
For a very interesting history of formula~\eqref{General
-Lord Brouncker}, see Berndt~[8, p. 145]. In addition, Lord Brouncker's continued fraction formula also plays an important role in Landau's constants, see~\cite{CXY,Cao1}.

The main aim in this subsection is to illustrate~(without proof) that the formula ~\eqref{General
-Lord Brouncker} is the fastest possible by making use of the method formulated in Sec.~4. Replacing $x$ by $4x+1$ in~\eqref{General
-Lord Brouncker} and then making some simple calculation, we obtain its equivalent forms as follows.
\begin{lem} Let $\Re x>-\frac 14$, we have
\begin{align}
\frac{\Gamma^2(x+\frac 12)}{\Gamma^2(x+1)}=\frac{1}{x+\frac 14+\frac{\frac{1}{32}}{x+\frac 14+\K_{m=1}^{\infty}\left(\frac{\frac{(2m+1)^2}{64}}{x+\frac 14}\right)}}.\label{General
-Lord Brouncker-1}
\end{align}
\end{lem}
\proof Now, we are in a position to treat the above formula directly. Let
\begin{align}
f(x)=\frac{\Gamma^2(x+\frac 12)}{\Gamma^2(x+1)}.
\end{align}
By the recurrence relation $\Gamma(x+1)=x \Gamma(x)$, we have
\begin{align}
\frac{f(x)}{f(x+1)}= \frac{(x+1)^2}{(x+\frac 12)^2}.\label{f-ratio-1}
\end{align}

By the Stirling's formula, it is not difficult to prove
\begin{align}
\lim_{x\rightarrow\infty}x^{b-a}\frac{\Gamma(x+a)}{\Gamma(x+b)}
=1.\label{Stirling's approximation}
\end{align}
Also see [1, p. 257, Eq. 6.1.47] or [8, p. 71, Lemma 2].

It follows readily from~\eqref{Stirling's approximation} that
\begin{align}
\lim_{x\rightarrow +\infty}x f(x)=1,\label{Brouncker-Initial-Con}
\end{align}
i.e, we take $\nu=1$ in~\eqref{Condition-1}.
\bigskip

\noindent {\bf(Step 1) The initial-correction.} According to~\eqref{Brouncker-Initial-Con}, we take $\Phi_0(x)=x+a$ for some constant $a$, to be specified below. From \eqref{f-ratio-1} and \eqref{MT-Mathematica-Program}, it is not difficult to prove that
\begin{align}
&\ln \frac{f(x)}{f(x+1)}+\ln\frac{\Phi_0(x)}{\Phi_0(x+1)}
=2\ln \frac{x+1}{x+\frac 12}+\ln\frac{x+a}{x+1+a}\\
=&\frac{-1/4 + a}{x^2}+O\left(\frac{1}{x^3}\right).\nonumber
\end{align}
Solve the equation $-1/4 + a=0$, we get $a=1/4$. By~\emph{Mortici-transformation}, we obtain
\begin{align}
\Phi_0(x)=x+\frac 14,\quad CF_0(f(x))=\frac{1}{x+\frac 14}.
\end{align}
As we need to use ~\emph{Mortici-transformation} in each correction-process, so will not mention it for the sake of simplicity.

\bigskip
\noindent {\bf(Step 2) The first-correction.} Let us expand the following function into a power series in terms of $1/x$:
\begin{align}
&\ln \frac{f(x)}{f(x+1)}+\ln\frac{\Phi_0(x)+\frac{\kappa_0}{x}}
{\Phi_0(x+1)+\frac{\kappa_0}{x+1}}
=2\ln \frac{x+1}{x+\frac 12}+\ln\frac{x+\frac 14+\frac{\kappa_0}{x}}
{x+\frac 54+\frac{\kappa_0}{x+1}}\\
=&\frac{-1/16 + 2\kappa_0}{x^3}+O\left(\frac{1}{x^4}\right).\nonumber
\end{align}
We solve the equation $-1/16 + 2\kappa_0=0$, and obtain $\kappa_0=1/32\neq 0$. Hence we take the first-correction $\mathrm{MC}_1(x)$ to be \emph{Type-I}, i.e.
\begin{align}
\mathrm{MC}_1(x)=\frac{\kappa_0}{x+\lambda_0}.
\end{align}
Since
\begin{align}
\ln \frac{f(x)}{f(x+1)}+\ln\frac{\Phi_0(x)+\frac{\kappa_0}{x+\lambda_0}}
{\Phi_0(x+1)+\frac{\kappa_0}{x+1+\lambda_0}}
=\frac{\frac{3}{128} - \frac{3 \lambda_0}{32}}{x^4}+O\left(\frac{1}{x^5}\right),\nonumber
\end{align}
we enforce $\frac{3}{128} - \frac{3 \lambda_0}{32}=0$, and deduce
$\lambda_0=\frac 14$. Thus,
\begin{align}
\mathrm{MC}_1(x)=\frac{\frac{1}{32}}{x+\frac 14},
\quad CF_1(f(x))=\frac{1}{x+\frac 14+\frac{\frac{1}{32}}{x+\frac 14}}.
\end{align}
\bigskip

\noindent {\bf(Step 3) The second-correction to the sixth-correction.} Now we take $\mathrm{MC}_2(x)$ to be \emph{Type-I}, and let
\begin{align}
\mathrm{MC}_2(x)=\frac {\kappa_0}{ x+\lambda_0+
\frac {\kappa_1}{x+\lambda_1}}.
\end{align}
By using\eqref{MT-Mathematica-Program}, we have
\begin{align}
&\ln \frac{f(x)}{f(x+1)}+\ln\frac{\Phi_0(x)+\mathrm{MC}_2(x)}
{\Phi_0(x+1)+\mathrm{MC}_2(x+1)}\\
=&\frac{\frac{9}{512} - \frac{\kappa_1}{8}}{x^5}+\frac{5 (-27 + 176 \kappa_1 + 64 \kappa_1 \lambda_1)}{2048 x^6}+O\left(\frac{1}{x^7}\right).\nonumber
\end{align}
Solve the equations
\begin{align}
\begin{cases}
\frac{9}{512} - \frac{\kappa_1}{8}=0,\\
-27 + 176 \kappa_1 + 64 \kappa_1 \lambda_1=0,
\end{cases}
\end{align}
we attain
\begin{align}
\kappa_1=\frac{9}{64},\quad \lambda_1=\frac 14.
\end{align}
We take the $k$th-correction $\mathrm{MC}_k(x)$ to be \emph{Type-I}, then repeat the above approach like the second-correction, and solve successively the coefficients $\kappa_j$ and $\lambda_j$ ($2\le j\le 6$) as follows:
\begin{align}
& \kappa_2=\frac{25}{64},\quad \lambda_2=\frac 14;\quad \kappa_3=\frac{49}{64},\quad \lambda_3=\frac 14;\quad
\kappa_4=\frac{81}{64},\quad \lambda_4=\frac 14;\\ &\kappa_5=\frac{121}{64},\quad \lambda_5=\frac 14;\quad \kappa_6=\frac{169}{64},\quad \lambda_6=\frac 14;\quad\kappa_7=\frac{225}{64},\quad \lambda_7=\frac 14.
\end{align}
From these results, it is not difficult to guess that
\begin{align}
\kappa_m=\frac{(2m+1)^2}{64},\quad \lambda_m=\frac 14.
\end{align}
Further, we apply~\eqref{MT-Mathematica-Program} to check that the above conjecture holds true for some larger $m$. In this way, we finally test that the fastest possible formula should be~\eqref{General
-Lord Brouncker-1}.\qed
%At same time, this approach provides a method for guessing
%a continued fraction expansion.

\subsection{The continued fraction formulas involving the volume of the unit ball}
Let $\Omega_n$ be defined by~\eqref{Omega-n-Def}. The main purpose of this subsection is to present the following two theorems.
\begin{thm} Let $n\ge 1$ be a positive integer. Then
\begin{align}
\frac{\Omega_{n}^2}{\Omega_{n-1}\Omega_{n+1}}=
\frac{2(n+1)}{2n+1+\K_{m=0}^{\infty}\left(
\frac{(2m+1)^2}{2(2n+1)}\right)}.
\end{align}
\end{thm}
\proof. It follows from~\eqref{Omega-n-Def} and the recurrence relation $\Gamma(x+1)=x \Gamma(x)$ that
\begin{align}
\frac{\Omega_{n}^2}{\Omega_{n-1}\Omega_{n+1}}=
\frac{\Gamma(\frac n2+\frac 12)\Gamma(\frac n2+\frac 32)}{\Gamma^2(\frac n2+1)}=\frac{n+1}{2}\frac{\Gamma^2(\frac n2+\frac 12)}{\Gamma^2(\frac n2+1)}.
\end{align}
Replacing $x$ by $\frac n2$ in~\eqref{General
-Lord Brouncker-1}, then after simplification, we get easily the desired assertion.\qed

\begin{thm} Let $n\in \mathbb{N}$, then
\begin{align}
\frac{\Omega_{n-1}}{\Omega_n}=\frac{1}{2\sqrt{\pi}}
\sqrt{2n+1+\K_{m=0}^{\infty}
\left(\frac{(2m+1)^2}{2(2n+1)}\right)
}.
\end{align}
\end{thm}
\proof From~\eqref{Omega-n-Def}, we have
\begin{align}
\frac{\Omega_{n-1}}{\Omega_n}=\frac{1}{\sqrt{\pi}}\frac{\Gamma(\frac n2+1)}{\Gamma(\frac n2+\frac 12)}.
\end{align}
Replacing $x$ by $\frac n2$ in~\eqref{General
-Lord Brouncker-1}, then taking reciprocals of both sides, finally substituting it into the above formula, this will complete the proof of Theorem 4.\qed

\begin{rem}
Condition~\eqref{Condition-1} is not an essential restriction. Actually, we can extend our method to any negative integer $\nu$.
For example, by taking reciprocals of both sides in~\eqref{General
-Lord Brouncker-1}, we have
\begin{align}
\frac{\Gamma^2(x+1)}{\Gamma^2(x+\frac 12)}=x+\frac 14+\frac{\frac{1}{32}}{x+\frac 14+}\K_{m=1}^{\infty}\left(\frac{\frac{(2m+1)^2}{64}}{x+\frac 14}\right),\quad \Re x>0.
\end{align}
In this case, we take $\nu=-1$. It should be remarked that we can discover the above formula directly by using an approach similarly to Lemma 5.
\end{rem}

\begin{rem}
To the best of our knowledge, formula~\eqref{General
-Lord Brouncker} and~\eqref{General
-Lord Brouncker-1} were
possibly neglected by many mathematicians for about more than twenty years, until 2013 I. Gavrea and M. Ivan mentioned it in their paper~\cite{GI}.
\end{rem}

\section{A continued fraction formula of Ramanujan}
The following lemma is Entry 39 in Berndt~[8, p. 159], which is one of three principal formulas involving gamma functions given by Ramanujan. It is very difficult for us to imagine how Ramanujan discovered those beautiful continued fraction formulas. Maybe our method provides a theoretical basis.
\begin{lem} Let $l$ and $n$ denote arbitrary complex numbers. Suppose that $x$ is complex with $\Re x >0$ or that either $n$ or  $l$ is an odd integer. Then
\begin{align}
P:=&\frac{\Gamma\left(\frac 14(x+l+n+1)\right)\Gamma\left(\frac 14(x-l+n+1)\right)\Gamma\left(\frac 14(x+l-n+1)\right)\Gamma\left(\frac 14(x-l-n+1)\right)}{\Gamma\left(\frac 14(x+l+n+3)\right)\Gamma\left(\frac 14(x-l+n+3)\right)\Gamma\left(\frac 14(x+l-n+3)\right)\Gamma\left(\frac 14(x-l-n+3)\right)}\label{Entry-39}\\
=&
\begin{array}{ccccccccccc}
8 & & 1^2-n^2 & & 1^2-l^2 & & 3^2-n^2 & & 3^2-l^2 &\\
\cline{1-1}\cline{3-3}\cline{5-5}\cline{7-7}\cline{9-9}\cline{11-11}
(x^2-l^2+n^2-1)/2 & + & 1 & + & x^2-1 & + & 1 & +&x^2-1 & +\cdots
\end{array}.\nonumber
\end{align}
\end{lem}

By replacing $x$ by $4x$, and taking $(l,n)=(0,0)$, $(l,n)=(1/4,1/2)$,
$(l,n)=(1/3,1/2)$, $(l,n)=(1/8,1/2)$, respectively, the authors have checked that Lemma 6 is not optimal continued fraction expansion. Now, by employing  these test, we may refine it in a uniform expression as follows.
\begin{thm} Under the same conditions of Lemma 6, we have
\begin{align}
P=&\begin{array}{ccccccc}
8& & (1^2-n^2)(1^2-l^2) & & (3^2-n^2)(3^2-l^2) &  \\
\cline{1-1}\cline{3-3}\cline{5-5}\cline{7-7}
\frac 12(x^2-l^2-n^2+1)&-&x^2-l^2-n^2+(3^2+1^2-1)&-&x^2-l^2-n^2+(5^2+3^2-1)
&-\cdots
\end{array}\\
=&\begin{array}{ccc}
8&\\
\cline{1-1}\cline{3-3}
\frac 12(x^2-l^2-n^2+1)&+
\end{array}
\K_{m=1}^{\infty}
\left(\frac{-\left((2m-1)^2-n^2\right)\left((2m-1)^2-l^2\right)}
{x^2-l^2-n^2+8 m^2+1}\right).\nonumber
\end{align}
\end{thm}
\proof
We follow the method of Entry 25 in Berndt~[8, p. 141]. First, we rewrite Lemma 6 in the form
\begin{align*}
\frac{8}{P}+\frac{1}{2}(x^2+l^2-n^2-1)=x^2-1+
\begin{array}{ccccccccc}
1^2-n^2 & & 1^2-l^2 & & 3^2-n^2 & & 3^2-l^2 &\\
\cline{1-1}\cline{3-3}\cline{5-5}\cline{7-7}\cline{9-9}
1 & + & x^2-1 & + & 1 & +&x^2-1 & +\cdots
\end{array},
\end{align*}
or
\begin{align*}
\frac{1}{8/P+\frac{1}{2}(x^2+l^2-n^2-1)}=
\begin{array}{ccccccccccc}
1&&1^2-n^2 & & 1^2-l^2 & & 3^2-n^2 & & 3^2-l^2 &\\
\cline{1-1}\cline{3-3}\cline{5-5}\cline{7-7}\cline{9-9}\cline{11-11}
x^2-1&+&1 & + & x^2-1 & + & 1 & +&x^2-1 & +\cdots
\end{array}.
\end{align*}
Secondly, by Entry 14 of Berndt~[8, p. 121]~(an infinity form see [8, p. 157]), we have
\begin{align}
&\frac{1}{8/P+\frac{1}{2}(x^2+l^2-n^2-1)}\\
=&
\begin{array}{ccccccc}
1& & (1^2-n^2)(1^2-l^2)& & (3^2-n^2)(3^2-l^2) &\\
\cline{1-1}\cline{3-3}\cline{5-5}\cline{7-7}
x^2-n^2&-&x^2-l^2-n^2+(3^2+1^2-1)&-&x^2-l^2-n^2+(5^2+3^2-1)&
-\cdots
\end{array}\nonumber\\
&
\begin{array}{ccc}
\quad&\left( (2m-1)^2-n^2\right)\left( (2m-1)^2-l^2\right)& \\
\cline{2-2}
-&x^2-l^2-n^2+\left( (2m+1)^2+(2m-1)^2-1\right)&-\cdots
\end{array}.\nonumber
\end{align}

Note that $(2m+1)^2+(2m-1)^2-1=8 m^2+1$. Now take the reciprocal of both sides above and then solve for $P$,
which again involves taking reciprocals. This will finish the proof of Theorem 5.\qed

The following theorem is the fastest possible form for Entry 26 in Berndt~[8, p. 145].
\begin{thm}
Suppose that either $n$ is an odd integer and $x$ is any complex number or that $n$ is an arbitrary complex number and $\Re x>0$. Then
\begin{align}
&\frac{\Gamma^2\left(\frac 14(x+n+1)\right)\Gamma^2\left(\frac 14(x-n+1)\right)}{\Gamma^2\left(\frac 14(x+n+3)\right)\Gamma^2\left(\frac 14(x-n+3)\right)}\\
=&\begin{array}{ccc}
8&\\
\cline{1-1}\cline{3-3}
\frac 12(x^2-n^2+1)&+
\end{array}
\K_{m=1}^{\infty}\left(
\frac{-(2m-1)^2\left((2m-1)^2-n^2\right)}
{x^2-n^2+8 m^2+1}\right).\nonumber
\end{align}
\end{thm}
\proof Set $l=0$ in Theorem 5, the desired equality follows at once.\qed

Similarly, we give another form of the Corollary in Berndt~[8, p. 146].

\begin{cor} If $\Re x>0$, then
\begin{align}
\frac{\Gamma^4\left(\frac 14(x+1)\right)}{\Gamma^4\left(\frac 14(x+3)\right)}
=\begin{array}{ccc}
8&\\
\cline{1-1}\cline{3-3}
\frac 12(x^2+1)&+
\end{array}
\K_{m=1}^{\infty}
\left(\frac{-(2m-1)^4}
{x^2+8 m^2+1}\right).
\end{align}
\end{cor}
\proof We set $n=0$ in Theorem 6, this completes the proof of the corollary readily.\qed

\section{Some new conjectural continued fraction formulas}
In this section, we will give three examples to illustrate how to guess their fastest possible continued fraction expansions. For the recent results involving these functions, see Mortici, Cristea and Lu~\cite{MCL}, Cao and Wang~\cite{CW}, and Chen~\cite{Chen}.

\subsection{For $\frac{\Gamma^3(x+\frac 13)}{\Gamma^3(x+1)}$}

In this subsection, we will use the function $\frac{\Gamma^3(x+\frac 13)}{\Gamma^3(x+1)}$ as an example to explain how to guess its fastest possible continued fraction expansion, which consists of the following steps.
\bigskip

{\bf (1).} Define
\begin{align*}
f(x):=\frac{\Gamma^3(x+\frac 13)}{\Gamma^3(x+1)}.
\end{align*}
Find the structure of $CF_k\left(f(x)\right)$ or $\mathrm{MC_k(x)}$ by \emph{Mortici transformation} and~\eqref{MT-Mathematica-Program}. We may determine that $CF_k\left(f(x)\right)$ has the form of \emph{Type-II}, and its $\mathrm{MC}$-point $\omega$ equals to $1/6$. Here we omit the details for finding those coefficients in $CF_k\left(f(x)\right)$, since the proof is very similar to that of Sec.~5 or Subsection 8.3 below.
\bigskip

{\bf (2).}  We denote $CF_k\left(f(x)\right)$ in the \emph{simplified form} like~\eqref{canonical-form}:
\begin{align}
CF_k\left(f(x)\right)=\frac{1}{(x+\omega)^2+\lambda_{-1}+
\K_{j=0}^{k-1}\left(\frac{\kappa_j}{(x+\omega)^2+\lambda_{j}}
\right)},
\end{align}
where $\lambda_{-1}=\frac{5}{2^2 3^3}$.
\bigskip

{\bf (3).}  We write two sequences $(\kappa_m)_{m\ge 0}$ and $(\lambda_m)_{m\ge 0}$ in the \emph{canonical form}, then extract their common factors, respectively. For example, one may use \emph{Mathematica} command ``FactorInteger" to do that. In this way, we denote these two
sequences in the form
\begin{align*}
&(\kappa_0,\kappa_1,\kappa_2,\kappa_3,\ldots)=-\frac{2}{3^6}
\left(\frac{1^3 1^3}{1^3},\frac{2^3 5^3}{3^2},\frac{4^3 7^3}{5^2},\frac{5^3 11^3}{7^2},
\frac{7^3 13^3}{9^2},\ldots\right),\\
&\lambda_{-1}=\frac{5}{2^2 3^3}, (\lambda_0,\lambda_1,\lambda_2,\lambda_3,\ldots)=\frac{1}{2^23^3}(
\frac{5^27}{1\cdot 3},\frac{3307}{3\cdot 5},\frac{17167}{5\cdot 7},\frac{5\cdot 31\cdot 353}{7\cdot 9},\ldots ).
\end{align*}

\bigskip

{\bf (4).}  Now we will look for the general terms of the sequences $(\kappa_m)_{m\ge 0}$ and $(\lambda_m)_{m\ge 0}$. We try to decompose them into some more simpler ``partial sequences".
\bigskip

(4-1). We observe easily that the sequence $(a_m)_{m\ge 0}=(1,3,5,\ldots)$ has the general term $a_m=2m+1$. While, for the sequence $(b_m)_{m\ge 0}=(1\cdot 3,3\cdot 5,5\cdot 7,\ldots)$, its general term is $b_m=(2m+1)(2m+3)$.
\bigskip

(4-2). Let us consider the sequence $(\alpha_m)_{m\ge 0}=(1,2,4,5,7,8,10,11,\ldots)$, which is the sequence generated by deleting the sequence $(3k)_{k\ge 1}$ from the positive integer sequence $(k)_{k\ge 1}$. We can check that the sequence $(\alpha_m)_{m\ge 0}$ satisfies the following difference equation
\begin{align*}
\alpha_m-\alpha_{m-1}=\begin{cases}
1,&\mbox{if $m$ is an odd,}\\
2,& \mbox{if $m$ is an even,}
\end{cases}
\end{align*}
with the initial value $\alpha_0=1$. Hence, we deduce that the general term equals to
\begin{align*}
\alpha_m=m+1+\lfloor\frac{m}{2}\rfloor.
\end{align*}
\bigskip

(4-3). Similarly, the sequence $(\beta_m)_{m\ge 0}=(1,5,7,11,13,17,19,\ldots)$ satisfies the following difference equation
\begin{align*}
\beta_m-\beta_{m-1}=\begin{cases}
4,&\mbox{if $m$ is an odd,}\\
2,& \mbox{if $m$ is an even,}
\end{cases}
\end{align*}
with the initial condition $\beta_0=1$, and we get
\begin{align*}
\beta_m=3 m+1+\frac{1-(-1)^m}{2}.
\end{align*}
\bigskip

(4-4). The sequence $(\xi_m)_{m\ge 0}=(5^27,3307,17167,5\cdot 31\cdot 353,5\cdot 13\cdot 2063,5\cdot 7\cdot 19\cdot 419,516847,7\cdot 13\cdot 9697, \ldots)$ is most difficult. Consider a new sequence
$(u_m)_{m\ge 0}$ to be defined by
\begin{align}
u_m:=\xi_m \mod (2m+1)(2m+3).
\end{align}
By using \emph{Mathematica} command ``mod", we can verify
\begin{align}
(u_0,u_1,u_2,\dots)=&(58, 220, 490, 868, 1354, 1948, 2650, 3460,\ldots)\\
=&2(29, 110, 245, 434, 677, 974, 1325, 1730,\ldots)\nonumber\\
:=&2 (v_m)_{m\ge 0}.\nonumber
\end{align}
\bigskip
We may check that the sequence
$(v_m)_{m\ge 0}$ satisfies the following difference equation
\begin{align}
v_{m}-2 v_{m-1}+v_{m-2}=108
\end{align}
with the initial conditions $v_0=29$ and $v_1=110$. Solve this difference equation of \emph{order} 2, we can deduce that the general term equals to
\begin{align}
v_m=27(m+1)^2+2.
\end{align}
Now we rewrite the general term $\xi_m$ in the form
\begin{align}
\xi_m=2(2m+1)(2m+3)v_m+w_m.
\end{align}
We may check that $(w_0,w_1,w_2,\ldots)=(1, 7, 17, 31, 49, 71, 97, 127,\dots)$. Quite similarly to the previous sequence $(v_m)_{m\ge 0}$, $(w_m)_{m\ge 0}$ also satisfies a difference equation as follows
\begin{align}
w_{m}-2 w_{m-1}+w_{m-2}=4
\end{align}
with the initial conditions $w_0=1$ and $w_1=7$. In this way, we get
\begin{align}
w_m=2 (m+1)^2-1.
\end{align}

Substituting (8.5) and (8.8) into (8.6), we discover
\begin{align}
\xi_m=2(2m+1)(2m+3)\left(27(m+1)^2+2\right)+2 (m+1)^2-1.
\end{align}
Combining the above results and after some simplification, we  conjecture that the general terms should be
\begin{align}
\kappa_m=&-\frac{2}{729}\frac{\left(m+1+\lfloor\frac{m}{2}
\rfloor\right)^3\left(3 m+1+\frac{1-(-1)^m}{2}
\right)^3}{(2m+1)^2},\quad (m\ge 0)\label{kappa-def}\\
\lambda_m=&\frac{1}{108}\left(2(27 (m+1)^2+2)+\frac{2 (m+1)^2-1}{(2m+1)(2m+3)}\right),\quad (m\ge -1).\label{lambda-def}
\end{align}
Note that we used the fact that the last formula also holds true for $m=-1$.
\bigskip

{\bf (5).} Define two sequences $(\kappa_m)_{m\ge 0}$ and $(\lambda_m)_{m\ge -1}$ by~\eqref{kappa-def} and~\eqref{lambda-def}, respectively. By making use of~\eqref{MT-Mathematica-Program}, we check that the above conjectures are still true for some ``larger" $m$.
\bigskip

{\bf (6).} Further simplification for the general term $\kappa_m$. Actually, we have
\begin{align}
\left(m+1+\lfloor\frac{m}{2}
\rfloor\right)\left(3 m+1+\frac{1-(-1)^m}{2}
\right)=&\frac{9}{8}\left((2 m + 1)^2 - (\frac 13)^2\right)\\
=&\frac{(3m+1)(3m+2))}{2},\nonumber
\end{align}
which may be proved easily according to $m$ is an odd and an even, respectively. Hence
\begin{align}
\kappa_m=-\frac{1}{2916}\frac{(3m+1)^3(3m+2)^3}{(2m+1)^2},\quad (m\ge 0).\label{kappa-def-1}
\end{align}
\bigskip

Finally, we propose the following reasonable conjecture.
\bigskip

\noindent{\bf Open Problem 1.} Let two sequences $(\kappa_m)_{m\ge 0}$ and $(\lambda_m)_{m\ge -1}$ be define by~\eqref{kappa-def-1} and~\eqref{lambda-def}, respectively. Let real $x> -1/6$, then we have
\begin{align}
\frac{\Gamma^3(x+\frac 13)}{\Gamma^3(x+1)}=\frac{1}{(x+\frac 16)^2+\lambda_{-1}+\frac{\kappa_0}{(x+\frac 16)^2+\lambda_0+\frac{\kappa_1}{(x+\frac 16)^2+\lambda_1+\frac{\kappa_2}{(x+\frac 16)^2+\lambda_2+\frac{\kappa_3}{(x+\frac 16)^2+\lambda_3+\ddots}}}}}.\label{Open Problem 1}
\end{align}
\begin{rem}
Open Problem 1 means that if there exists a fastest possible continued fraction expansion for the function $\frac{\Gamma^3(x+\frac 13)}{\Gamma^3(x+1)}$, then it must be the continued fraction expression of the right side in~\eqref{Open Problem 1}.
\end{rem}
\bigskip

Replacing $x$ by $x-1/6$ and then after some simplification, we get the following equivalent forms of Open Problem 1.
\bigskip

\noindent{\bf Open Problem $1^{\prime}$.} Let real $x> 0$, then
\begin{align}
\frac{\Gamma^3(x+\frac 16)}{\Gamma^3(x+\frac 56)}
=\begin{array}{ccc}
1&\\
\cline{1-1}\cline{3-3}
x^2+\frac{5}{108}&+
\end{array}\K_{n=1}^{\infty}\left(\frac{-\frac{(3n-2)^3(3n-1)^3}
{2916(2n-1)^2}}
{x^2+\frac{1}{108}\left(2(27 n^2+2)+\frac{2 n^2-1}{(2n-1)(2n+1)}\right)}\right).
\end{align}

\subsection{For $\frac{\Gamma^3(x+\frac 23)}{\Gamma^3(x+1)}$}
The main purpose of this subsection is to conjecture  the fastest possible continued fraction expansion for the function $f(x)$, which is defined by
 $$f(x):=\frac{\Gamma^3(x+\frac 23)}{\Gamma^3(x+1)}.$$
We follow the same method described in last subsection. By testing, we observe that $CF_k\left(f(x)\right)$ has the form of \emph{Type-I}, and its $\mathrm{MC}$-point $\omega$ is $1/3$. Some computation data
are listed as follows:
\begin{align}
CF_k\left(f(x)\right)=\frac{1}{x+\frac 13+\K_{j=0}^{k-1}\left(\frac{\kappa_j}{x+\frac 13}\right)},
\end{align}
where
\begin{align}
\kappa_0=\frac{1}{27},\quad (\kappa_1,\kappa_2,\kappa_3,\ldots)=\frac{1}{54}\left(\frac{2^3}{1},
\frac{4^3}{3},\frac{5^3}{3},\frac{7^3}{5},
\frac{8^3}{5},\frac{10^3}{7},\frac{11^3}{7},
\frac{13^3}{9},\frac{14^3}{9},\frac{16^3}{11},\frac{17^3}{11},
\ldots\right).\label{kappa-1}
\end{align}
\bigskip
Similarly to the sequence $(\alpha_m)_{m\ge 0}$ in last subsection, it is not difficult to verify that the general term of
the sequence $(\lambda_m)_{m\ge 1}$ should be
\begin{align}
\lambda_m=\frac{1}{54}\frac{\left(m+1+\lfloor\frac{m}{2}
\rfloor\right)^3}{2\lfloor\frac{m}{2}\rfloor+1}.\label{kappa-2}
\end{align}

\noindent{\bf Open Problem 2.} For all real $x>-\frac 13$, we have
\begin{align}
\frac{\Gamma^3(x+\frac 23)}{\Gamma^3(x+1)}=\frac{1}{x+\frac 13 +\frac{\frac{1}{27}}{x+\frac 13+\frac{\frac{1}{54}\frac{2^3}{1}}{x+\frac 13+\frac{\frac{1}{54}\frac{4^3}{3}}{x+\frac 13+\frac{\frac{1}{54}\frac{5^3}{3}}{x+\frac 13+\ddots}}}}}.
\end{align}
Replace $x$ by $x-1/3$, we have the following equivalent forms of Open Problem 2.
\bigskip

\noindent{\bf Open Problem $2^{\prime}$.} Let $\kappa_0=\frac {1}{27}$, and the sequence $(\lambda_m)_{m\ge 1}$ be defined as~\eqref{kappa-2}. Let $x>0$, then
\begin{align}
\frac{\Gamma^3(x+\frac 13)}{\Gamma^3(x+\frac 23)}=\begin{array}{ccc}
1&\\
\cline{1-1}\cline{3-3}
x&+
\end{array}
\K_{m=0}^{\infty}\left(\frac{\kappa_m}{x}\right).\label{Open Problem-2-1}
\end{align}
Since the partial coefficients of the continued fraction of the right side in~\eqref{Open Problem-2-1} are all positive, we can prove the following consequence easily.
\begin{cor} Let $x>0$. Assume that Open Problem $2^{\prime}$ is true, then for all non-negative integer $k$
\begin{align}
\begin{array}{ccc}
1&\\
\cline{1-1}\cline{3-3}
x&+
\end{array}
\K_{j=0}^{2k+1}\left(\frac{\kappa_j}{x}\right)
<\frac{\Gamma^3(x+\frac 13)}{\Gamma^3(x+\frac 23)}<\begin{array}{ccc}
1&\\
\cline{1-1}\cline{3-3}
x&+
\end{array}
\K_{j=0}^{2k}\left(\frac{\kappa_j}{x}\right).
\end{align}
\end{cor}
\begin{rem}
The authors have checked that Corollary 3 is true for $k\le 10$.
\end{rem}

\subsection{For $\frac{\Gamma(x+\eta)\Gamma(x+1-\eta)}{\Gamma^2(x+1)}$}
Let $\eta$ be a real number with $0<\eta<1$. In this subsection, we will
discuss the continued fraction approximation for the ratio of the gamma functions
\begin{align}
G_{\eta}(x):=\frac{\Gamma(x+\eta)\Gamma(x+1-\eta)}{\Gamma^2(x+1)}.
\label{G-eta-Def}
\end{align}

It follows from~\eqref{Stirling's approximation} that
\begin{align}
\lim_{x\rightarrow\infty}x G_{\eta}(x)=1.
\end{align}
Now let us begin to look for $CF_k(G_{\eta}(x))$.
\bigskip

\noindent {\bf(Step 1) The initial-correction.} Note that $\nu=1$ in~\eqref{Mortici-transformation}. It follows readily from the recurrence formula $\Gamma(x+1)=x\Gamma(x)$ that
\begin{align*}
\frac{G_{\eta}(x)}{G_{\eta}(x+1)} =\frac{(x+1)^2}{(x+\eta)(x+1-\eta)}.
\end{align*}

Now we apply \emph{Mortici-transformation} to determine $\Phi_0(x)=x+c_0$. By making use of \emph{Mathematica} software, one has
\begin{align}
\ln\frac{(x+1)^2}{(x+\eta)(x+1-\eta)}+\ln \frac{x+c_0}{x+1+c_0}
=\frac{c_0-\eta+ \eta^2}{x^2}+O\left(\frac{1}{x^3}\right).
\end{align}
Solve $c_0-\eta+ \eta^2=0$, we obtain $c_0=\eta - \eta^2$.
\bigskip

\noindent {\bf(Step 2) The first-correction.} Let
\begin{align}
\mathrm{MC}_1(x)=\frac{\kappa_0}{x+\lambda_0},
\end{align}
similarly to the initial-correction, we also have
\begin{align}
&\ln\frac{(x+1)^2}{(x+\eta)(x+1-\eta)}+\ln \frac{x+c_0+\mathrm{MC}_1(x)}{x+1+c_0+\mathrm{MC}_1(x+1)}\\
=&\frac{2 \kappa_0 - \eta^2 + 2 \eta^3 - \eta^4}{x^3}+\frac{-3 \kappa_0 - 3 \kappa_0 \lambda_0 - 3 \kappa_0 \eta + 2 \eta^2 + 3 \kappa_0 \eta^2 -
 3 \eta^3 - \eta^4 + 3 \eta^5 - \eta^6}{x^4}\nonumber\\
&+\frac{g(x)}{x^5}
+O\left(\frac{1}{x^6}\right),\nonumber
\end{align}
where
\begin{align*}
g(x)=&4 \kappa_0 - 2 \kappa_0^2 + 6 \kappa_0 \lambda_0 + 4 \kappa_0 \lambda_0^2 + 6 \kappa_0 \eta + 4 \kappa_0 \lambda_0 \eta - 3 \eta^2 - 2 \kappa_0 \eta^2 \\
&- 4 \kappa_0 \lambda_0 \eta^2 + 4 \eta^3 - 8 \kappa_0 \eta^3 + 2 \eta^4 + 4 \kappa_0 \eta^4 - 2 \eta^5 - 4 \eta^6 + 4 \eta^7 - \eta^8.
\end{align*}
Solve
\begin{align}
\begin{cases}
&2 \kappa_0 - \eta^2 + 2 \eta^3 - \eta^4=0\\
&-3 \kappa_0 - 3 \kappa_0 \lambda_0 - 3 \kappa_0 \eta + 2 \eta^2 + 3 \kappa_0 \eta^2 -
 3 \eta^3 - \eta^4 + 3 \eta^5 - \eta^6=0\\
&g(x)=0,
\end{cases}
\end{align}
we get
\begin{align}
\kappa_0 = \frac{(-1 + \eta)^2 \eta^2}{2},\qquad \lambda_0 =
\frac{1 - \eta + \eta^2}{3}.
\end{align}
\bigskip

\noindent {\bf(Step 3) The second-correction to the sixth-correction.} Similarly to the first-correction, we use \emph{Mathematica} software to
find that the second-correction to the sixth-correction are the form of \emph{Type-I}, and then solve all coefficients in these correction functions. It should be remarked that there is a parametric $\eta$, the related computations will become very huge and complex. So we need to manipulate \emph{Mathematica} command ``{\emph{Simplify} }" . Here we list the final computing results as follows:
\begin{align*}
&\kappa_1 =\frac{(-2 - \eta + \eta^2)^2}{36}
,\qquad \lambda_1 =\frac{4 - \eta + \eta^2}{15}
;\\
& \kappa_2 =\frac{(-6 - \eta + \eta^2)^2}{100}
,\qquad \lambda_2 =\frac{9 - \eta + \eta^2}{35}
;\\
& \kappa_3 =\frac{(-12 - \eta + \eta^2)^2}{196}
,\qquad \lambda_3 =\frac{16 - \eta + \eta^2}{63}
;\\
& \kappa_4 =\frac{(-20 - \eta + \eta^2)^2}{324}
,\qquad \lambda_4 =\frac{25 - \eta + \eta^2}{99};\\
& \kappa_5 =\frac{(-30 - \eta + \eta^2)^2}{484}
,\qquad \lambda_5 =\frac{36 - \eta + \eta^2}{143}.
\end{align*}

Similarly to Open Problem 1, by careful data analysis and further checking, we may propose the following conjecture.
\bigskip

\noindent {\bf Open Problem 3.} For $x>0$, then
\begin{align}
\frac{\Gamma(x+\eta)\Gamma(x+1-\eta)}{\Gamma^2(x+1)}=\frac{1}
{x+\eta-\eta^2+\K_{m=0}^{\infty}\left(\frac{\kappa_m}{x+\lambda_m}
\right)},
\end{align}
where $\kappa_0=\frac{(-\eta+\eta^2)^2}{2}$ and
\begin{align}
\kappa_m=&\frac{\left(-m (m + 1) - \eta+ \eta^2\right)^2}{4 (2 m + 1)^2}=\frac{(m+\eta)^2(m-\eta+1)^2}{4 (2 m + 1)^2},\quad m\ge 1\\
\lambda_m=&\frac{(1 + m)^2 - \eta + \eta^2}{(2 m + 1) (2 m + 3)},\quad m\ge 0.
\end{align}
\begin{rem}
If we take $\eta=\frac 12$, then the above conjecture implies \eqref{General
-Lord Brouncker-1}~(i.e. the generalized Lord Brouncker's continued fraction formula). Here we note that
\begin{align*}
&\frac{(m+\eta)^2(m-\eta+1)^2}{4 (2 m + 1)^2}=\frac{(m+\frac 12)^4}{4 (2 m + 1)^2}=\frac{(2 m + 1)^2}{2^6},\\
&\frac{(1 + m)^2 - \eta + \eta^2}{(2 m + 1) (2 m + 3)}=
\frac{(1 + m)^2 - (\frac 12)^2}{(2 m + 1) (2 m + 3)}=\frac{(m+\frac 12)(m+\frac 32)}{(2 m + 1) (2 m + 3)}=\frac 14.
\end{align*}
\end{rem}

\section{Conclusions}
In this paper, we present a systematical way to construct a best possible finite and infinite continued fraction approximations for a class of functions. In particular, the method  described in Sec.~4 is suitable for the ratio of the gamma functions, e.g. many examples can be found in the nice survey papers Qi~\cite{Qi} and Qi and Luo~\cite{QL}. As our method is constructive, so all involving computations may be manipulated by a suitable symbolic computation software, e.g.~\emph{Mathematica}. In some sense, the main advantage of our method is that such formal continued fraction approximation of order $k$ is the fastest possible when $x$ tends to infinity. Concerning applications in approximation theory, numerical computation, our method represents a much better approximation formula than the power series approach~(e.g. Taylor's formula) for a kind of ``good functions".

In addition, the \emph{multiple-correction method} provides a useful tool for testing and guessing the continued fraction expansion involving a specified function. So our method should help advance the approximation theory, the theory of continued fraction, the generalized hypergeometric function, etc. Further, if we can obtain some new continued fraction expansions, probably these formulas could be used to study the irrationality, transcendence of the involved constants.

\end{document}